\documentclass[11pt]{article}
\pdfoutput=1
\usepackage[margin=1in]{geometry}
\usepackage{mathptmx}
\usepackage{titlesec}
\usepackage{etoolbox}
\usepackage{setspace}
\usepackage{graphicx}
\usepackage{natbib}
\usepackage{newtxtext}
\usepackage{newtxmath}
\usepackage{hyperref}
\usepackage{float}
\usepackage{soul}
\hypersetup{
	colorlinks = true,
	urlcolor   = blue,
	citecolor  = black,
}

\newtheorem{definition}{Definition}
\newtheorem{theorem}{Theorem}

\newtheorem{proof}{Proof}

\newcommand{\RomanNumeralCaps}[1]
\linenumbers

\title{Optimization Based Data Enrichment Using Stochastic Dynamical System 
Models}

\begin{document}
	\maketitle
	Griffin M. 
	Kearney$\mathrm{^{a,b}}$\footnote{griffin.kearney@opbdatainsights.com}, 
	Makan Fardad$\mathrm{^{b}}$\footnote{makan@syr.edu}
	
	$\mathrm{^a}${\it OpB Data Insights LLC, Syracuse, NY 13224}
	
	$\mathrm{^b}${\it Syracuse University, Syracuse, NY 13244}
	
	\vskip\bigskipamount 
	\vskip\medskipamount
	\leaders\vrule width \textwidth\vskip0.4pt 
	\vskip\bigskipamount 
	\vskip\medskipamount
	\nointerlineskip
	\begin{abstract}
		We develop a general framework for state estimation in systems 
		modeled with noise-polluted continuous time dynamics and discrete 
		time noisy measurements.
		Our approach is based on maximum likelihood estimation and employs
		the calculus of variations to derive optimality conditions for 
		continuous time functions.
		We make no prior assumptions on the form of the mapping from 
		measurements to state-estimate or on the distributions of the noise
		terms, making the framework more general than Kalman filtering/smoothing
		where this mapping is assumed to be linear and the noises Gaussian.
		The optimal solution that arises is interpreted as a continuous time
		spline, the structure and temporal dependency of which is determined
		by the system dynamics and the distributions of the process and 
		measurement noise.
		Similar to Kalman smoothing, the optimal spline yields increased data
		accuracy at instants when measurements are taken, in addition to 
		providing continuous time estimates outside the measurement instances.
		We demonstrate the utility and generality of our approach via 
		illustrative examples that render both linear and nonlinear data
		filters depending on the particular system.
	\end{abstract}
	
	\section{Introduction}
	\label{sec:intro}
	In most practical applications, measurements are noisy observables sampled 
	from 
	underlying dynamics that are governed by physical laws and polluted by 
	stochastic processes. Knowledge of both the governing temporal dynamics and 
	the 
	distribution of the uncertainty can be exploited to add accuracy to the 
	measurements and counter the effect of noise. It also allows for 
	interpolation 
	and extrapolation, to obtain estimates for those time instances at which no 
	measurements were made. Here we collectively refer to data accuracy 
	improvements and completion as data enrichment.
	
	In this work, a novel technique for the enrichment of data governed 
	by stochastic dynamics and observed through noisy measurement systems is 
	developed.
	We adopt the term data enrichment to describe the mitigation 
	of errors in measurements (I.E. the reduction of noise) and the estimation 
	of 
	missing data points (I.E. interpolation / extrapolation of measurements) 
	simultaneously.
	For instance, if the position of a vehicle were periodically measured in 
	time 
	using a noisy GPS system, then data enrichment would improve the accuracy 
	of 
	the GPS measurements at the sampling times and provide vehicle 
	position estimates for the times in which no GPS samples have been taken.
	
	Our new approach achieves this goal through joint analysis of stochastic 
	dynamical systems and stochastic measurements.
	Data accuracy improvements and completion are achieved by generating 
	maximum 
	likelihood estimates using the distributions that characterize the 
	stochastic  
	dynamics and measurements.
	In the present work, since our attention is on dynamical systems, our method
	results in optimal splines over time. 
	
	A spline refers to a piecewise function with a prescribed structure on 
	each of 
	the underlying sub-intervals that support it over the domain.
	Spline techniques are well-established and among them cubic spline 
	interpolation is used across a variety of 
	applications \citep{Durrleman1989,Hsieh1978,Gesemann2018}.
	However, the classical construction of splines often requires subjective 
	assumptions to be enforced in the fitting process which lack physical 
	justification.
	In our treatment this is not the case. 
	The spline shape and the constraints enforced to optimally enrich the data 
	are 
	determined entirely by the stochastic dynamics and measurements within the 
	maximum likelihood estimation framework.
	In other words, our technique automatically provides optimally structured 
	splines for each 
	system of interest.
	Moreover, we demonstrate that cubic splines are the optimal spline shape 
	for a simple, common class of stochastic dynamics, 
	thus providing new novel insights that support their use.
	
	Kalman filtering is an alternative well-established method for data 
	enrichment.
	There are similarities between the reasoning behind our technique and the 
	widely used Kalman approach 
	\citep{Kalman1960,Kalman1961,Barratt2020,Gustafsson2011,Kulikov2018,Lewis2017}.
	Namely, Kalman methods provide estimates by balancing dynamic and 
	measurement stochastic process likelihoods simultaneously.
	Despite this high-level similarity, the present optimization method and the 
	Kalman method are fundamentally different.
	Kalman assumes, from the outset, that the filter acts linearly on 
	observed measurements whereas our method requires no such restriction.
	This limitation results in sub-optimal filter performance outside of 
	systems 
	involving 
	linear dynamical equations and Gaussian processes 
	\citep{Hendeby2005}.
	Indeed, the optimality of the linear filter follows naturally from our 
	framework 
	when considering these systems without such an assumption; it is an 
	automatic 
	consequence of our technique.
	
	We proceed with development of the new method by first describing required 
	preliminaries in Section \ref{sec:model-overview-and-prelims}.
	During the preliminary discussion we motivate and establish a critical 
	extension of instantaneous process distributions to continuous intervals 
	which 
	is required for forming the general problem.
	In Section \ref{sec:motivation-example} we provide a simple motivating 
	example 
	to show how the problem is formed by considering the motion of a
	point mass, but refrain from computing the solution at this stage.
	The example is provided to make the general treatment of the problem in 
	Section 
	\ref{sec:MLSR} more approachable, and it is within 
	this discussion that we develop the optimization conditions in detail.
	Finally, in Section \ref{sec:applications} we examine illustrative examples 
	that demonstrate how optimal splines are constructed using the main theorem 
	that results from the general analysis.
	
	A summary of the contributions of this work are as follows:
	\begin{enumerate}
		\item We develop a novel optimization framework to generate continuous 
		time 
		enrichments of discretely sampled systems with known stochastic 
		dynamical process and measurement model distributions. 
		This development includes a new method for extending instantaneous 
		process 
		distributions to continuous intervals.
		\item We establish the necessary optimization criteria which 
		characterize 
		solutions of the new framework. 
		The optimization criteria are provided as a set of 
		general ordinary differential equations which determine the optimal 
		spline 
		structure, and a set of algebraic equations which characterize the 
		boundary 
		conditions at measurement points.
		These conditions are derived through application of the Calculus of 
		Variations. 
		\item When the governing dynamics are linear, and both the dynamical 
		process and measurements are subject to Gaussian stochastic processes, 
		then we show that the new theory results in optimal splines that behave 
		as 
		linear filters on noisy measurements. 
		We demonstrate that this produces cubic splines in a simple 
		illustrative 
		set of dynamics (simple particle motion), 
		but that it also produces non-cubic splines in a second illustrative 
		case 
		(simple harmonic motion).
		We then analyze an example on the general analysis of the optimal 
		splines for Kalman systems (linear dynamics with Gaussian stochastic 
		distributions) which concretely demonstrates that the linear optimal 
		solution flows naturally from the new theory without prior assumption. 
		Finally, to demonstrate generality, we conclude with consideration of 
		two examples that produce non-linear optimal solutions. 
		The first of these examines an example of non-Gaussian process noise, 
		and the second involves non-linear governing dynamics.
	\end{enumerate}
	
	\section{Model Overview and Preliminaries}
	\label{sec:model-overview-and-prelims}
	The variable $x$ is used to represent the state of a system, and the 
	variable 
	$y$ 
	is used to represent measurements observed within it.%
	\footnote{Throughout this work when a vector-valued function of time,
		say $x$, appears without explicit reference to $t$ then it implies the 
		function
		over its entire temporal support, i.e. $x(\cdot)$.}
	We begin by assuming a general form of governing dynamics on the state $x$,
	\begin{equation}
		\label{eqn:general-dynamics}
		\dot{x}(t) = f(t,x(t)) + \nu(t,v(t)),
	\end{equation}
	and a general form of measurements
	\begin{equation}
		\label{eqn:general-measurements}
		y(t) = g(t,\dot{x}(t)) + h(t,x(t)) + \xi(t,w(t)).
	\end{equation}
	For each $t$ the vector $x(t)$ is assumed to be of dimension $n_x$ and the
	vector $y(t)$ of dimension $n_y$.
	The functions $f$, $g$, $h$, $\nu$, and $\xi$ are permitted to be 
	arbitrary, 
	emphasizing that in this framework we allow for non-linear dynamics and 
	measurements.
	Moreover, these functions are also allowed to have explicit dependence on 
	time.
	
	In general, $v$ and $w$ are used to denote stochastic processes, they are 
	vectors of 
	dimensions $n_v$ and $n_w$ respectively, and we assume 
	them to be statistically independent of one another.
	At each instant $t$, $v(t)$ and $w(t)$ denote random variables, and both 
	are 
	associated with corresponding probability density functions 
	$\rho_v(t,v(t))$ 
	and $\rho_w(t,w(t))$, which characterize their distributions.
	Furthermore we assume that the images of $v$ and $w$ are statistically 
	independent of themselves at distinct time instants.
	In other words, $v(t_0)$ is independent of $v(t_1)$ for $t_1 \neq t_0$.
	
	The dynamical stochastic process (\ref{eqn:general-dynamics}) is assumed to 
	hold for all times $t$. 
	The stochastic measurement description (\ref{eqn:general-measurements}) is 
	only defined 
	for a strict and possibly continuous subset of time instants.
	We denote the entire time horizon
	under consideration as $\mathscr{T}$, and 
	the subset where measurements have been captured as $\mathscr{T}_M 
	\subseteq \mathscr{T}$.
	
	In many practical problems $\mathscr{T}_M$ will be a finite collection of 
	disjoint discrete points, but continuous measurements on intervals are also 
	permissible.
	Our goal will be to quantify the likelihood of outcomes on the set 
	$\mathscr{T}$, and for the points in $\mathscr{T}_M$ this is determined by 
	both 
	the dynamical process and measurement randomness as modeled through $v$ and 
	$w$.
	On the other hand, for the set $\mathscr{T} - \mathscr{T}_M$ the likelihood 
	is 
	only determined by $v$ as no measurements are taken at these times.
	We have assumed that the dynamics and measurement stochastic processes are 
	independent, which allows us to write the piecewise function
	\begin{equation}
		\label{eqn:distribution-breakdown}
		\rho(t) = 
		\begin{cases}
			\rho_v(t,v(t)) & t \in \mathscr{T} - \mathscr{T}_M\\
			\rho_v(t,v(t)) \rho_w(t,w(t)) & t \in \mathscr{T}_M,
		\end{cases}
	\end{equation}
	for any time $t \in \mathscr{T}$.
	The function $\rho$ assigns a real number, the 
	probability 
	density, to the value $v(t)$ (or to the pair $v(t)$, $w(t)$ when 
	measurements 
	are present) at the time $t$.
	It will be important to keep in mind, once we form the maximum
	likelihood problem, that $\rho$ is a function of $v$ and $w$ and depends on 
	their probability distributions even though these dependencies are not made 
	explicit in our choice of notation.
	
	Computing probabilities at discrete points is possible using this 
	function directly.
	Suppose that $t_0$ and $t_1$ are two time instants of interest at which no 
	measurements 
	are taken.
	We assume that the process values $v(t)$ are independent for distinct 
	times, 
	and therefore the joint distribution of $v(t_0)$ and $v(t_1)$ is expressed 
	as
	\begin{align*}
		\rho({\{t_0,t_1\}}) & = \rho(t_0) \rho(t_1)\\
		& = \rho_v(t_0,v(t_0)) \rho_v(t_1,v(t_1)),  
	\end{align*}
	a simple product.
	This can be done for a set of more than two discrete points by including 
	more 
	factors in the product, and adjusted using 
	(\ref{eqn:distribution-breakdown}) 
	for times 
	when measurements are taken.
	
	Difficulties arise applying this approach when we must compute the 
	probability 
	over intervals instead of discrete points.
	Assume that $\tau = (t_0,t_1)$ is an open interval over which no 
	measurements 
	are 
	taken.
	It is tempting to write
	\[
	\rho(\tau) = \prod_{t \in (t_0,t_1)} \rho_v(t,v(t)),
	\]
	but it is unclear how the product in this expression 
	should be computed or if it is even well-defined as it consists of an 
	uncountably infinite number of factors.
	
	We resolve this difficulty by proposing a method for extending 
	instantaneous 
	distributions in the Appendix (Section 
	\ref{subsec:extend-instant-dynamics}).
	The treatment provides analytic machinery for extending functions described 
	by 
	$\rho$ to functionals on continuous 
	intervals that we denote with $\mu(\tau,\rho)$.
	The extension is of the form
	\begin{equation}
		\label{eqn:mu-computation}
		\mu(\tau,\rho) = e^{\frac{1}{|\tau|} \int_{\tau} \mathrm{ln} \rho(t) 
			\hspace{0.10cm} dt},
	\end{equation}
	where $|\tau|$ denotes the length of the interval,
	and as a point of convention if $\tau = \{t\}$, then we 
	define
	\[
	\mu(\tau,\rho) = \rho(t)
	\]
	for simplicity of notation.
	
	The extension $\mu(\tau,\rho)$ assigns a real number to the image of the 
	process $v$ (or to the pair of images $v$, $w$ on intervals with 
	measurements) on $\tau$.
	We interpret $\mu$ as an analog of the probability density, but on 
	interval/process pairs 
	instead of discrete time/random variable pairs.
	In alignment with this interpretation, we define how the joint likelihood 
	over 
	two separate intervals is treated to maintain the parallel.
	We say that two intervals $\tau_0$ and $\tau_1$ are separated if there 
	exists 
	some point between them.
	For instance $(a,b)$ and $(b,c)$ are separated for all $a < b < c$ since 
	the 
	point $\{b\}$ lies between them, but $(a,b]$ 
	and $(b,c)$ are not separated.
	If $\tau_0$ and $\tau_1$ are separated intervals, then we require the 
	extension 
	function satisfies the relation
	\begin{equation}
		\label{eqn:mu-separated}
		\mu(\tau_0 \cup \tau_1, \rho) = \mu(\tau_0,\rho) \mu(\tau_1, \rho).
	\end{equation}
	
	In all problems considered here, the sets $\mathscr{T}_M$ and 
	$\mathscr{T}-\mathscr{T}_M$ will permit representations as unions of 
	separated 
	intervals and discrete points.
	The maximum likelihood problem we propose has an objective of the general 
	form
	\begin{equation}
		\label{eqn:gen-MLE-ob}
		J = \mu(\mathscr{T}-\mathscr{T}_M,\rho) \mu(\mathscr{T}_M,\rho),
	\end{equation}
	which is maximized subject to the governing dynamics 
	(\ref{eqn:general-dynamics})--(\ref{eqn:general-measurements}).
	Equation (\ref{eqn:mu-separated}) is used to simplify each of the two terms 
	on the right-hand-side of
	(\ref{eqn:gen-MLE-ob}) for a given problem.
	We recall from (\ref{eqn:distribution-breakdown}) that $\rho$ depends on 
	$v$, 
	$w$ and their probability distributions.
	Loosely speaking, the maximum likelihood problem seeks to find the most 
	likely 
	trajectories for $v$ and $w$ that satisfy the constraints 
	(\ref{eqn:general-dynamics})--(\ref{eqn:general-measurements}). 
	
	In Section \ref{sec:MLSR} we fully develop and analyze the general problem. 
	Before we do this however, we use the next section to formulate the 
	optimization 
	problem on an illustrative example to motivate the general treatment.
	
	\section{Motivational Example}
	\label{sec:motivation-example}
	We focus attention on a simple scenario.
	We consider the one dimensional undamped motion of a point mass under 
	Gaussian stochastic forcing. 
	Position measurements are collected periodically 
	in time during the particle's motion, which are subject to Gaussian errors.
	
	The governing dynamics of the system are
	\begin{equation}
		\label{eqn:motivation-dynamics}
		\ddot{r}(t) = v(t)
	\end{equation}
	where $r(t)$ is a real number representing the position at time $t$ and 
	$v$ is a scalar 
	stochastic process with
	\begin{equation}
		\label{eqn:motivation-rho-v}
		\rho_v(t,v(t)) = \frac{1}{\sqrt{2 \pi} \sigma_p} e^{-\frac{v(t)^2}{2 
				\sigma_p^2}}.
	\end{equation}
	We denote the process noise variance as $\sigma_p^2$ in 
	(\ref{eqn:motivation-rho-v}) and assume that it is known.
	The position measurements are modeled
	\begin{equation}
		\label{eqn:motivation-meas}
		y(t) = r(t) + w(t)
	\end{equation}
	where $y(t)$ denotes the noisy measurement at time $t$ and $w$ is a 
	stochastic 
	process with
	\begin{equation}
		\label{eqn:motivation-rho-w}
		\rho_w(t,w(t)) = \frac{1}{\sqrt{2 \pi} \sigma_m} e^{-\frac{w(t)^2}{2 
				\sigma_m^2}}.
	\end{equation}
	We denote the measurement noise variance as $\sigma_m^2$ and also assume 
	that 
	it is known.
	The state space vector is written as
	\[
	x(t) = \begin{pmatrix}
		r(t)\\
		\dot{r}(t)
	\end{pmatrix}
	\]
	and enables us to express this system in standard form:
	\begin{align}
		\label{eqn:motivation-state}
		\dot{x}(t) & = 
		\begin{pmatrix}
			0 & 1\\
			0 & 0
		\end{pmatrix} x(t) + 
		\begin{pmatrix}
			0\\1
		\end{pmatrix}v(t)\\
		\label{eqn:motivation-measurements}
		y(t) & = \begin{pmatrix}
			1 & 0
		\end{pmatrix} x(t) + w(t).
	\end{align}
	
	Suppose that measurements are taken periodically in time with sampling 
	frequency $f_0$.
	We represent the sampling times with the set $\{t_k\}_{k=0}^K$ and assume
	\[
	t_{k+1} - t_k = \frac{1}{f_0}
	\]
	for $k = 0,1,\ldots,K\!-\!1$.
	Let us consider the time interval $\mathscr{T} = [t_0,t_K]$, where 
	$t_0$ is the first and $t_K$ is the final measurement time.
	We have 
	\[
	\mathscr{T}_M = \{t_k\}_{k=0}^K,
	\] 
	by assumption, and
	\[
	\mathscr{T} - \mathscr{T}_M = \cup_{k=0}^{K-1} (t_k,t_{k+1})
	\]
	as a result.
	The set $\mathscr{T}_M$ consists of separated discrete points, and 
	$\mathscr{T}-\mathscr{T}_M$ consists of separated continuous intervals.
	This allows us to use relation (\ref{eqn:mu-separated}) to form the 
	objective (\ref{eqn:gen-MLE-ob}) for this example.
	Namely,
	\begin{equation}
		\label{eqn:motivation-dyn-likelihood}
		\mu(\mathscr{T}-\mathscr{T}_M,\rho) = \prod_{k=0}^{K-1} 
		\mu((t_k,t_{k+1}),\rho)
	\end{equation}
	and
	\begin{equation}
		\label{eqn:motivation-meas-likelihood}
		\mu(\mathscr{T}_M,\rho) = \prod_{k=0}^{K} \mu(\{t_k\},\rho).
	\end{equation}
	Equations (\ref{eqn:motivation-rho-v}) and (\ref{eqn:motivation-rho-w}) are 
	used in (\ref{eqn:distribution-breakdown}) to define $\rho$ in these 
	expressions.
	
	Maximizing the objective, computed using (\ref{eqn:gen-MLE-ob}), is 
	equivalent 
	to maximizing its logarithm.
	We reformulate the objective as
	\begin{equation}
		\label{eqn:motivation-obj}
		\mathrm{ln} J = \sum_{k=0}^{K-1} \mathrm{ln} 
		\mu((t_k,t_{k+1}),\rho) + \sum_{k=0}^K \mathrm{ln}\mu(\{t_k\},\rho)
	\end{equation}
	after substitution using (\ref{eqn:motivation-dyn-likelihood}) and 
	(\ref{eqn:motivation-meas-likelihood}) for practicality.
	Furthermore, a final simplification of (\ref{eqn:motivation-obj}) is 
	performed 
	using 
	(\ref{eqn:mu-computation}) to write
	\begin{equation}
		\label{eqn:motivation-obj-final}
		\mathrm{ln} J = f_0 \sum_{k=0}^{K-1} \int_{t_k}^{t_{k+1}} \mathrm{ln} 
		\rho(t) 
		dt + \sum_{k=0}^K \mathrm{ln} \rho(t_k),
	\end{equation}
	where $f_0$ is the sampling frequency of the measurements.
	
	The governing optimization problem is formed by maximizing 
	(\ref{eqn:motivation-obj-final}) 
	subject to constraints which enforce the state and measurement models in 
	(\ref{eqn:motivation-state}) and (\ref{eqn:motivation-measurements}).
	This is equivalent to the problem
	\begin{align}
		\notag
		\underset{x,v,w}{\mathrm{min.}} & \hspace{0.1cm} f_0 \sum_{k=0}^{K-1}  
		\int_{t_k}^{t_{k+1}}  
		\frac{v(t)^2}{2 \sigma_p^2} 
		dt + \sum_{k=0}^K (\frac{v(t_k)^2}{2 \sigma_p^2}+\frac{w(t_k)^2}{2 
			\sigma_m^2})\\
		\label{eqn:motivation-optimization-prob}
		\mathrm{s.t.} & \hspace{0.1cm} \dot{x}(t) = 
		\begin{pmatrix}
			0 & 1\\
			0 & 0
		\end{pmatrix} x(t) + 
		\begin{pmatrix}
			0\\1
		\end{pmatrix}v(t) \hspace{0.25cm} \text{ for } t \in \mathscr{T}\\
		\notag
		& \hspace{0.1cm} y(t) = \begin{pmatrix}
			1 & 0
		\end{pmatrix} x(t) + w(t) \hspace{0.25cm} \text{ for } t \in 
		\mathscr{T}_M,
	\end{align}
	in this simple example.
	We postpone the analysis of this problem, instead examining a general 
	treatment 
	in Section \ref{sec:MLSR} that is specifically applied to this system in 
	Section \ref{sec:applications}. 
	However, we preview the resulting spline which arises from its solution to 
	build intuition before moving forward.
	Interestingly, the optimal solution in this system takes the form of a 
	cubic 
	spline. 
	We demonstrate the optimal spline behavior on simulated data in the 
	following subsection, and discuss additional simulations of this system 
	in Section \ref{sec:cubic-splines-as-MLSR}.
	\subsubsection*{Preview (Cubic Splines)}
	\label{subsubsec:cubic-splines}
	
	We simulate the point mass dynamics with noisy measurements.
	We set $\sigma_p = 4$ and $\sigma_m = 1$ for the purposes of this 
	demonstration.
	The data consists of a simulated trajectory, constructed using discrete 
	updates 
	such that
	\[
	r(t + dt) = r(t) + \dot{r}(t) dt + \frac{1}{2} a(t) dt^2,
	\]
	and
	\[
	\dot{r}(t + dt) = \dot{r}(t) + a(t) dt.
	\]
	For each simulation time $t$, the simulated acceleration $a(t)$ is randomly 
	drawn from a zero-mean Gaussian distribution with variance $\sigma_p^2$, 
	and 
	$dt$ was set to $0.01$ time units.
	
	The trajectory is generated using a starting time $t_0 = 0$ and a final 
	time 
	$t_K = 10$, with
	$r(0) = 10$ and $\dot{r}(0) = 0$ used to initialize the trajectory.
	Simulated measurements are collected with sampling frequency $f_0 = 5 
	\frac{\text{samples}}{\text{time unit}}$ using the 
	relation
	\[
	y(t_k) = r(t_k) + w(t_k).
	\]
	The measurement noise $w(t_k)$ is randomly drawn from a zero-mean 
	Gaussian distribution with variance $\sigma_m^2$ for each measurement.
	
\begin{figure}
		\centering
		\includegraphics[scale=0.65]{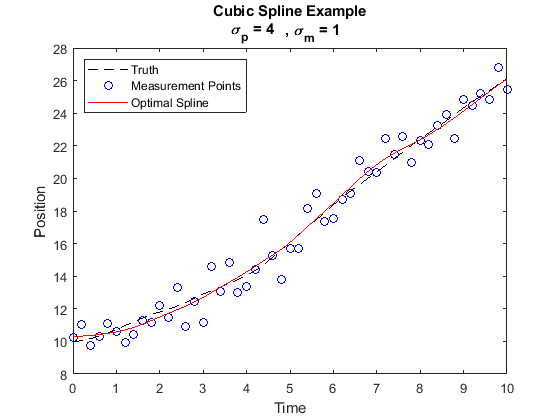}
		\caption{A simulated trajectory with noisy measurements and an optimal 
			cubic spline solution.}
		\label{fig:poly-spline-example}
	\end{figure}

	The simulated noisy measurements are used to construct an optimal spline 
	by solving the optimization problem (\ref{eqn:motivation-obj-final}).
	This spline is shown to be a cubic spline in the solution process which is 
	described later in Section \ref{sec:cubic-splines-as-MLSR}. 
	The noiseless trajectory, the simulated measurement points, and the optimal 
	cubic spline are all shown in Figure \ref{fig:poly-spline-example}.
	Notice that the optimization does not force the spline to pass through 
	each of the measurement points.
	A common constraint in spline applications is to require that the spline 
	agree 
	with measurements at sampling times.
	However, in the presence of substantial measurement noise this requirement 
	can 
	induce oscillations which are not present in the underlying true 
	trajectory, 
	causing accuracy breakdown.
	On the other hand, the present optimization theory provides both the cubic 
	structure without prior assumption and stitching conditions from each 
	interval 
	to the next automatically.
	The derivation of the cubic structure and optimal stitching conditions are 
	developed in detail in Section \ref{sec:cubic-splines-as-MLSR}. 
	This produces a spline that remains close to the true trajectory despite 
	the 
	measurement noise.
	
	\section{General Problem and Analysis}
	\label{sec:MLSR}
	In the general treatment we continue to restrict our attention to 
	$\mathscr{T}$ and  
	$\mathscr{T}_M$ as in 
	the motivational example.
	That is to say, we assume measurements are taken at discrete times with 
	some sampling frequency $f_0$ and we focus on enriching the data between 
	the 
	first and final measurements, inclusively.
	However, we modify the concepts utilized to formulate
	(\ref{eqn:motivation-optimization-prob}) to allow general
	stochastic distributions, dynamics, and measurements. 
	The resulting problem is written as
	\begin{align}
		\notag
		\underset{x,v,w}{\mathrm{max.}} & \hspace{0.1cm} f_0 \sum_{k=0}^{K-1}  
		\int_{t_k}^{t_{k+1}} 
		\mathrm{ln} 
		\rho(t) 
		dt + \sum_{k=0}^K \mathrm{ln} \rho({t_k})\\
		\label{eqn:general-optimization-prob}
		\mathrm{s.t.} & \hspace{0.1cm} \dot{x}(t) = f(t,x(t)) + \nu(t,v(t)) 
		\text{ 
			for } t \in \mathscr{T}\\
		\notag
		&  \hspace{0.1cm} y(t) = g(t,\dot{x}(t)) + h(t,x(t)) + \xi(t,w(t)) 
		\text{ for } t \in 
		\mathscr{T}_M.
	\end{align}
	
	Analysis of the problem proceeds in two steps.
	First, we transform (\ref{eqn:general-optimization-prob}) into an 
	unconstrained 
	optimization problem through use of
	the Langrangian and dual variables.
	Second, we apply techniques from the Calculus of Variations (CoV) to derive 
	necessary conditions on optimal solutions.
	Applying the CoV produces two types of optimal conditions: ordinary 
	differential equations that must be satisfied on $\mathscr{T} - 
	\mathscr{T}_M$, 
	and algebraic relations which must be satisfied on $\mathscr{T}_M$.
	The solution to the differential equations produces the spline structure 
	automatically, and the algebraic relations provide boundary conditions to 
	constrain the general solutions and in doing so fit the spline to a 
	specific 
	data set.
	
	It is useful to introduce compressed notation for clarity in development of 
	the 
	Lagrangian.
	We define
	\begin{equation}
		\label{eqn:general-compressed-dynamics}
		\phi(t) = \dot{x}(t) - f(t,x(t)) - \nu(t,v(t)),
	\end{equation}
	and
	\begin{equation}
		\label{eqn:general-compressed-measurements}
		\psi(t) = y(t) - g(t,\dot{x}(t)) - h(t,x(t)) - 
		\xi(t,w(t)).
	\end{equation}
	We introduce two sets of dual variables $\lambda$ and $\eta$.
	The first, $\lambda$, is a function of time defined on $\mathscr{T}$ and 
	relates 
	to the first set of constraints involving the dynamics
	(\ref{eqn:general-compressed-dynamics}).
	The second, $\eta$, is a discrete function of time defined only for times 
	in 
	$\mathscr{T}_M$ and corresponds to the second set of constraints relating 
	to 
	the measurements 
	(\ref{eqn:general-compressed-measurements}).
	Problem (\ref{eqn:general-optimization-prob}) is reformulated using the 
	Lagrangian objective
	\begin{align}
		\label{eqn:general-lagrange-problem}
		L = & f_0 \sum_{k=0}^{K-1} \int_{t_k}^{t_{k+1}} \mathrm{ln} 
		\rho(t) 
		- \lambda(t)^T \phi(t) dt 
		+ \sum_{k=0}^K \mathrm{ln} \rho({t_k}) - 
		\lambda(t_k)^T \phi(t_k) - \eta(t_k)^T \psi(t_k),
	\end{align}
	which yields the equivalent unconstrained optimization
	\begin{equation}
		\label{eqn:general-lagrangian-opt-problem}
		\underset{x,v,w,\lambda,\eta}{\mathrm{max.}} \hspace{0.1cm} L.
	\end{equation}
	
	The optimal conditions that are necessary for a valid solution contain many 
	Jacobian matrices and one-sided limits.
	We define notation conventions that allow us to express these concepts 
	efficiently.
	If $f$ is a differentiable vector-valued function from 
	$\mathbb{R}^n$ to $\mathbb{R}^m$, then we write
	\[
	f(u + \delta u) = f(u) + \frac{\partial f}{\partial u}(u) \delta u + 
	\mathrm{h.o.t}.
	\]
	
	In this expression $\frac{\partial f}{\partial u}$ is a linear operator of 
	dimension $m \times n$ and is the Jacobian of 
	$f$ with respect to $u$.
	If $x$ is a function of time, then we define for $\varepsilon > 0$
	\[
	x(t^-) = \lim_{\varepsilon \rightarrow 0} x(t - \varepsilon),
	\hspace{0.15in}
	x(t^+) = \lim_{\varepsilon \rightarrow 0} x(t + \varepsilon)
	\]
	as one-sided limits.
	
	The optimization conditions which characterize solutions of 
	(\ref{eqn:general-lagrangian-opt-problem}) are stated in Theorem 
	\ref{thm:Optimization-Conditions}.
	\begin{theorem}
		\label{thm:Optimization-Conditions}
		Assume $x$ is continuous for all $t \in \mathscr{T}$. If $x$, $v$, 
		$w$, 
		$\lambda$, and $\eta$ are an optimal solution of 
		(\ref{eqn:general-lagrangian-opt-problem}), then the following must 
		hold:
		\begin{align}
			\label{eqn:gen-opt-cond1}
			& \frac{\partial \mathrm{ln} \rho_v}{\partial v} + \frac{\partial 
				\nu}{\partial v}^T \lambda(t) = 0\\
			\label{eqn:gen-opt-cond2}
			& \dot{\lambda}(t) + \frac{\partial f}{\partial x}^T \lambda(t) = 
			0\\
			\label{eqn:gen-opt-cond3}
			& \dot{x}(t) - f(t,x(t)) - \nu(t,v(t)) = 0,
		\end{align}
		for all $t \in \mathscr{T} - \mathscr{T}_M$, and
		\begin{align}
			\label{eqn:gen-opt-cond4}
			&\frac{\partial \mathrm{ln} \rho_v}{\partial v} + \frac{\partial 
				\nu}{\partial v}^T \lambda(t_k) = 0\\
			\label{eqn:gen-opt-cond5}
			& \frac{\partial \mathrm{ln} \rho_w}{\partial w} + \frac{\partial 
				\xi}{\partial w}^T \eta(t_k) = 0\\
			\label{eqn:gen-opt-cond6}
			& \lambda(t_k) - \frac{\partial g}{\partial \dot{x}}^T \eta(t_k) = 
			0\\
			\label{eqn:gen-opt-cond7}
			& \frac{\partial f}{\partial x}^T \lambda(t_k) + \frac{\partial 
				h}{\partial x}^T 
			\eta(t_k) + f_0\lambda(t_k^+) - f_0\lambda(t_k^-) = 
			0\\
			\label{eqn:gen-opt-cond8}
			& y(t_k) - g(t_k,\dot{x}(t_k)) - h(t_k,x(t_k)) - \xi(t_k,w(t_k)) = 
			0\\
			\label{eqn:gen-opt-cond9}
			& \dot{x}(t_k) - f(t_k,x(t_k)) - \nu(t_k,v(t_k)) = 0,
		\end{align}
		for all $t_k \in \mathscr{T}_M$\footnote{When $\mathscr{T} = 
			[t_0,t_K]$, 
			then define $\lambda(t_0^-) = \lambda(t_K^+) = 
			0$ for simplicity}.
	\end{theorem}
	The proof of Theorem \ref{thm:Optimization-Conditions} is provided by
	standard application of the CoV.
	This generates a natural derivation of the full set of 
	optimization conditions listed within the theorem.
	An outline of the proof is provided in the Appendix 
	(Section \ref{subsec:Application of the Calculus of Variations});
	for technical details of the CoV the reader is referred to
	\citep{Clarke2013}.
	
	Equations (\ref{eqn:gen-opt-cond1})--(\ref{eqn:gen-opt-cond3}) define the 
	differential equations that are satisfied on each interval between 
	measurements.
	This system of equations depends only on $x$, $\lambda$, and $v$. 
	Therefore we conclude that the general structure of the optimal spline is 
	independent of the measurements, and it only depends on the dynamical 
	system 
	and the nature of the corresponding stochastic process.
	This is an interesting result, establishing that a system of stochastic 
	dynamics comes equipped with an optimal spline structure, prior to the 
	collection of measurements.
	
	Equations (\ref{eqn:gen-opt-cond4})--(\ref{eqn:gen-opt-cond9}) describe 
	the discrete conditions that must hold from one measurement-free interval 
	to 
	the next.
	This system of algebraic equations provides the constraints for the general 
	solutions to the differential equations, and its solution ultimately 
	enriches 
	the specific set of data by optimally fitting the spline to the 
	measurements.
	
	Constraint (\ref{eqn:gen-opt-cond7}) provides one-sided
	limits of $\lambda$ at each measurement point, and these values are used as 
	boundary conditions to fit solutions of 
	(\ref{eqn:gen-opt-cond2}).
	The dual variable $\lambda$ need not be continuous at the measurement 
	points in 
	general, and therefore the right and left handed limits are critical.
	On the other hand, under the assumption that $x$ is continuous at each of 
	the 
	measurement points, 
	we have
	\begin{equation}
		\label{eqn:gen-opt-cond-x-limit}
		x(t_k^\pm) = x(t_k)
	\end{equation}
	for all $t_k \in \mathscr{T}_M$.
	The values $\{x(t_k)\}_{k=0}^K$ are solved for using the algebraic 
	constraints, and 
	these values provide the boundary conditions 
	for (\ref{eqn:gen-opt-cond3}).
	The assumption that $x$ is continuous is a natural, physically motivated 
	one;
	if $f$ and $\nu$ are bounded then $x$ is necessarily continuous given 
	(\ref{eqn:general-dynamics}).
	
	\section{Illustrative Examples}
	\label{sec:applications}
	
	We provide illustrative examples, where the first two
	are chosen as specific instances of the third.
	We focus our attention in the first three examples on linear dynamical 
	systems with Gaussian stochastic processes.
	These are often good approximations of real systems,
	and are the standards assumptions under which Kalman filters/smoothers are 
	derived.
	These examples demonstrate that our theory produces optimal splines that 
	are 
	automatically adapted to the dynamics of an underlying system.
	
	For the first case, we return to the analysis of the motivational example 
	of Section \ref{sec:motivation-example}.
	We show that the optimal spline is cubic, develop 
	the fitting conditions, and exhibit its behavior on various simulated noisy 
	measurements.
	For the second case, we complicate the dynamics by considering a harmonic 
	oscillator and repeat the analysis.
	In this example, we show explicitly that the structure of the optimal 
	spline 
	radically changes with the dynamics.
	The piece-wise structure of the spline becomes a type of modified 
	harmonic, and is no longer cubic.
	Finally, for the third case we consider the Kalman smoothing
	problem for general linear systems.
	We compute the specific form of Theorem \ref{thm:Optimization-Conditions} 
	and show that the linear form of the filter, which is an 
	explicit assumption in Kalman
	filtering theory, is merely a natural consequence
	of our developed framework.
	
	We complete the illustrative examples by considering two additional 
	scenarios in which the optimal filter is non-linear.
	The first of these problems considers a case with a non-Gaussian process, 
	and the second includes a problem with non-linear dynamics.
	In these examples the linear Kalman approach is no longer optimal and 
	indeed a non-linear solution arises automatically from the new framework. 
	
	\subsection{Simple Particle Motion}
	\label{sec:cubic-splines-as-MLSR}
	We restate the dynamics (\ref{eqn:motivation-state}), measurements 
	(\ref{eqn:motivation-measurements}), and stochastic distributions
	(\ref{eqn:motivation-dyn-likelihood})--(\ref{eqn:motivation-meas-likelihood})
	here for convenience:
	\[
	\dot{x}(t) = 
	\begin{pmatrix}
		0 & 1\\
		0 & 0
	\end{pmatrix} x(t) + 
	\begin{pmatrix}
		0\\1
	\end{pmatrix}v(t),
	\]
	\[
	y(t) = \begin{pmatrix}
		1 & 0
	\end{pmatrix} x(t) + w(t),
	\]
	\[
	\rho_v(t,v(t)) = \frac{1}{\sqrt{2 \pi} \sigma_p} 
	e^{-\frac{1}{2 \sigma_p^2}v(t)^2},
	\]
	and
	\[
	\rho_w(t,w(t)) = \frac{1}{\sqrt{2 \pi} \sigma_m} 
	e^{-\frac{1}{2 \sigma_m^2}w(t)^2}.
	\]
	We apply Theorem \ref{thm:Optimization-Conditions} to derive the 
	optimal splines and stitching conditions.
	We begin with the differential equations describing the optimal spline 
	structure in (\ref{eqn:gen-opt-cond1})--(\ref{eqn:gen-opt-cond3}):
	\begin{align*}
		&-\frac{1}{\sigma_p^2} v(t) + \begin{pmatrix}
			0 & 1
		\end{pmatrix} \lambda(t) = 0\\
		& \dot{\lambda}(t) + \begin{pmatrix}
			0 & 0\\
			1 & 0\\
		\end{pmatrix} \lambda(t) = 0\\
		& \dot{x}(t) - \begin{pmatrix}
			0 & 1\\
			0 & 0
		\end{pmatrix} x(t) - \begin{pmatrix}
			0\\1
		\end{pmatrix}v(t) = 0.
	\end{align*}
	This must hold for all $t$ in intervals of the form $(t_k,t_{k+1})$.
	The solution on the interval $(t_k,t_{k+1})$ is of the general form
	\[
	v(t) = a_k (t - t_k) + b_k,
	\]
	\[
	\lambda(t) = \begin{pmatrix}
		a_k\\ a_k (t - t_k) + b_k
	\end{pmatrix},
	\]
	and
	\[
	x(t) = \begin{pmatrix}
		\frac{a_k \sigma_p^2}{6}(t - t_k)^3 + \frac{b_k \sigma_p^2}{2}(t - 
		t_k)^2 + c_k (t - t_k) + 
		d_k\\
		\frac{a_k \sigma_p^2}{2}(t - t_k)^2 + b_k \sigma_p^2 (t - t_k) + c_k
	\end{pmatrix},
	\]
	where $a_k$, $b_k$, $c_k$, and $d_k$ are a set of arbitrary constants.
	Notice that the particle position, the first component of $x(t)$, has a 
	cubic 
	form.
	This establishes explicitly how cubic splines arise naturally in the 
	stochastic 
	system defined by (\ref{eqn:motivation-state}) and 
	(\ref{eqn:motivation-dyn-likelihood}) in this framework.
	As an added benefit, the second component of $x(t)$ provides an optimal 
	estimate of the particle velocity.
	
	The algebraic conditions in 
	(\ref{eqn:gen-opt-cond4})--(\ref{eqn:gen-opt-cond9})
	are written as
	\begin{align*}
		& -\frac{1}{\sigma_p^2} v(t_k) + \begin{pmatrix}
			0 & 1
		\end{pmatrix} \lambda(t_k) = 0\\
		& -\frac{1}{\sigma_m^2} w(t_k) + \eta(t_k) = 0\\
		&\lambda(t_k) = 0\\
		&\begin{pmatrix}
			0 & 0\\1 & 0
		\end{pmatrix}\lambda(t_k) + 
		\begin{pmatrix}
			1\\0 
		\end{pmatrix}\eta(t_k) + f_0 \lambda(t_k^+) - f_0 \lambda(t_k^-) = 0\\
		&y(t_k) - \begin{pmatrix}
			1 & 0
		\end{pmatrix} x(t_k) - w(t_k) = 0\\
		&\dot{x}(t_k) - \begin{pmatrix}
			0 & 1\\0 & 0
		\end{pmatrix}x(t_k) - \begin{pmatrix}
			0\\1
		\end{pmatrix}v(t_k) = 0,
	\end{align*}
	for each $t_k \in \mathscr{T}_M$.
	This set of algebraic equations is used to constrain the general constants 
	on the intervals to construct the spline which best enriches the 
	measurements.
	A desirable attribute of splines is that they provide \emph{infinite 
		resolution} while only requiring \emph{finite representations}.
	Our optimal spline produces continuous time estimates but is simply 
	represented 
	by a finite set of $4K$ constants.
	
	\begin{figure}[h!]
		\centering
		\includegraphics[scale=0.65]{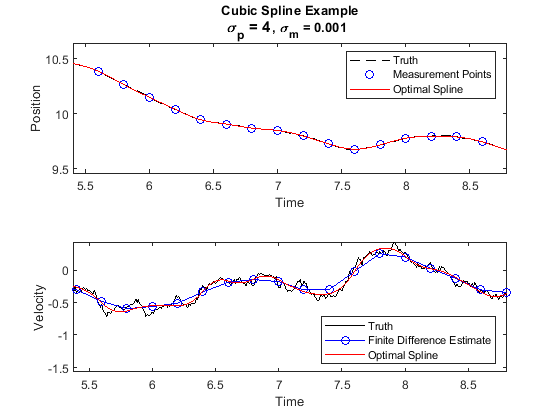}
		\caption{Simulated point mass trajectory with high accuracy 
			measurements. 
			[Top] Particle Position. [Bottom] Particle Velocity.}
		\label{fig:low-error-cubic-spline}
	\end{figure}
	\begin{figure}[h!]
		\centering
		\includegraphics[scale=0.65]{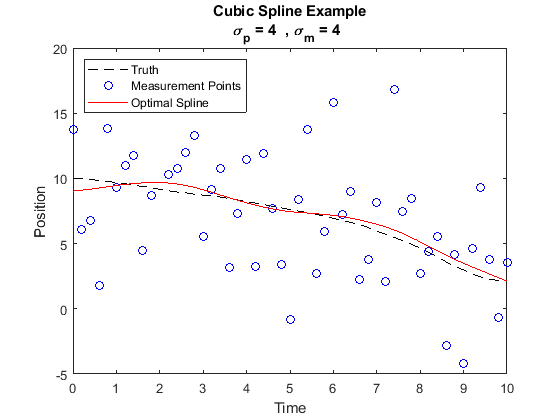}
		\caption{Simulated point mass trajectory with low accuracy 
			measurements.}
		\label{fig:high-error-cubic-spline}
	\end{figure}
	
	In Figure \ref{fig:low-error-cubic-spline} we show a simulated set of high 
	accuracy 
	measurements where $\sigma_m \ll \sigma_p$.
	The simulation was generated using the same procedure as in Section 
	\ref{subsubsec:cubic-splines} for the motivational example.
	The lower portion of Figure \ref{fig:low-error-cubic-spline} depicts the 
	true velocity, the velocity estimate arising from a piece-wise 
	linear fit (finite differences), and the velocity estimate automatically 
	produced by the new optimal spline.
	In Figure \ref{fig:high-error-cubic-spline} we show results from a 
	simulation 
	with much lower 
	measurement accuracy, i.e., with larger $\sigma_m$ than in Figure 
	\ref{fig:low-error-cubic-spline}.
	Compare the behavior of the optimal spline in this second simulation with 
	that 
	of the first.
	The optimal spline no longer adheres to each measurement point 
	individually, 
	but it does continue to approximate the true trajectory well.
	
	Importantly, the constants 
	$\{a_k,b_k,c_k,d_k\}_{k=0}^K$ are linearly related to the measured 
	observations.
	Therefore, the optimally estimated state at any time $t$ for 
	this system is linear with respect to the measurements.
	This is consistent with Kalman filtering theory, where it is known that a 
	linear filter 
	is the optimal choice for linear dynamical systems subject to Gaussian 
	processes.
	In our treatment the linear dependence of the estimated state on the 
	measured 
	values is a natural consequence following
	from the optimization framework rather than being imposed from
	the outset. 
	
	\subsection{Harmonic Oscillators}
	\label{subsec:harmonic-oscillator}
	We repeat the example in the previous section with the minor modification of
	replacing the governing dynamics by a harmonic oscillator with natural 
	frequency 
	$\omega^2$.
	The modified system is written as
	\begin{equation*}
		\label{eqn:harmonic-osc-state}
		\dot{x}(t) = 
		\begin{pmatrix}
			0 & 1\\
			-\omega^2 & 0
		\end{pmatrix} x(t) + 
		\begin{pmatrix}
			0\\1
		\end{pmatrix}v(t),
	\end{equation*}
	\begin{equation*}
		y(t) = \begin{pmatrix}
			1 & 0
		\end{pmatrix} x(t) + w(t),
	\end{equation*}
	\begin{equation*}
		\rho_v(t,v(t)) = \frac{1}{\sqrt{2 \pi} \sigma_p} 
		e^{-\frac{1}{2 \sigma_p^2}v(t)^2},
	\end{equation*}
	and
	\begin{equation*}
		\rho_w(t,w(t)) = \frac{1}{\sqrt{2 \pi} \sigma_m} 
		e^{-\frac{1}{2 \sigma_m^2}w(t)^2}.
	\end{equation*}
	
	Again we apply Theorem \ref{thm:Optimization-Conditions}. 
	First, examining the differential equations, we consider
	\begin{align*}
		&-\frac{1}{\sigma_p^2} v(t) + \begin{pmatrix}
			0 & 1
		\end{pmatrix} \lambda(t) = 0\\
		& \dot{\lambda}(t) + \begin{pmatrix}
			0 & -\omega^2\\
			1 & 0\\
		\end{pmatrix} \lambda(t) = 0\\
		& \dot{x}(t) - \begin{pmatrix}
			0 & 1\\
			-\omega^2 & 0
		\end{pmatrix} x(t) - \begin{pmatrix}
			0\\1
		\end{pmatrix}v(t) = 0.
	\end{align*}
	The solution on the interval $(t_k,t_{k+1})$ is given by
	\[
	v(t) = a_k \sin \omega t + b_k \cos \omega t,
	\]
	\[
	\lambda(t) = a_k \begin{pmatrix}
		-\omega \cos \omega t\\ \sin \omega t
	\end{pmatrix} + b_k \begin{pmatrix}
		\omega \sin \omega t\\ \cos \omega t
	\end{pmatrix},
	\]
	and
	\begin{align*}
		&x(t) = \begin{pmatrix}
			(\frac{\sigma_p^2 b_k}{2 \omega} t + c_k) \sin \omega t + 
			(-\frac{\sigma_p^2 a_k}{2 \omega}t + d_k) \cos \omega t\\
			(\frac{\sigma_p^2 a_k}{2}t + \frac{\sigma_p^2 b_k}{2 \omega} - 
			\omega 
			d_k)\sin \omega t + (\frac{\sigma_p^2 b_k}{2}t -\frac{\sigma_p^2 
				a_k}{2 
				\omega} + \omega c_k) \cos \omega t
		\end{pmatrix}
	\end{align*}
	where $a_k$, $b_k$, $c_k$, and $d_k$ are again arbitrary constants arising 
	in 
	the general solution.
	It is interesting to note the different structure of this solution compared 
	to the previous example, resulting from a small 
	modification of the dynamics.
	In particular, this optimal spline is a piece-wise function 
	of the form
	\[
	(\frac{\sigma_p^2 b_k}{2 \omega} t + c_k) \sin \omega t + 
	(-\frac{\sigma_p^2 a_k}{2 \omega}t + d_k) \cos \omega t,
	\]
	which we refer to as a \emph{modified harmonic}. 
	
\begin{figure}[h!]
		\centering
		\includegraphics[scale=0.65]{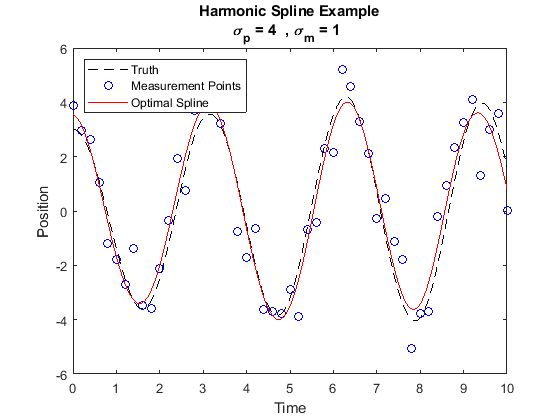}
		\caption{Simulated stochastic harmonic oscillator trajectory with 
			optimal 
			spline.}
		\label{fig:harmonic-spline}
	\end{figure}

	The algebraic conditions for this system are written as
	\begin{align*}
		& -\frac{1}{\sigma_p^2} v(t_k) + \begin{pmatrix}
			0 & 1
		\end{pmatrix} \lambda(t_k) = 0\\
		& -\frac{1}{\sigma_m^2} w(t_k) + \eta(t_k) = 0\\
		&\lambda(t_k) = 0\\
		&\begin{pmatrix}
			0 & -\omega^2\\1 & 0
		\end{pmatrix}\lambda(t_k) + 
		\begin{pmatrix}
			1\\0 
		\end{pmatrix}\eta(t_k) + f_0 \lambda(t_k^+) - f_0 \lambda(t_k^-) = 0\\
		&y(t_k) - \begin{pmatrix}
			1 & 0
		\end{pmatrix} x(t_k) - w(t_k) = 0\\
		&\dot{x}(t_k) - \begin{pmatrix}
			0 & 1\\-\omega^2 & 0
		\end{pmatrix}x(t_k) - \begin{pmatrix}
			0\\1
		\end{pmatrix}v(t_k) = 0,
	\end{align*}
	for each $t_k \in \mathscr{T}_M$.
	These equations are used to 
	constrain the general solution and construct a spline for the 
	specific data set.
	In this example we still have a system of linear dynamics and measurements 
	with Gaussian stochastic processes, and thus the new framework again 
	derives a linear relationship between the measured data and optimal spline.
	We simulate a stochastic harmonic oscillator and perform a representative 
	data enrichment in Figure \ref{fig:harmonic-spline}.
	
	\begin{figure}[h!]
		\centering
		\includegraphics[scale=0.55]{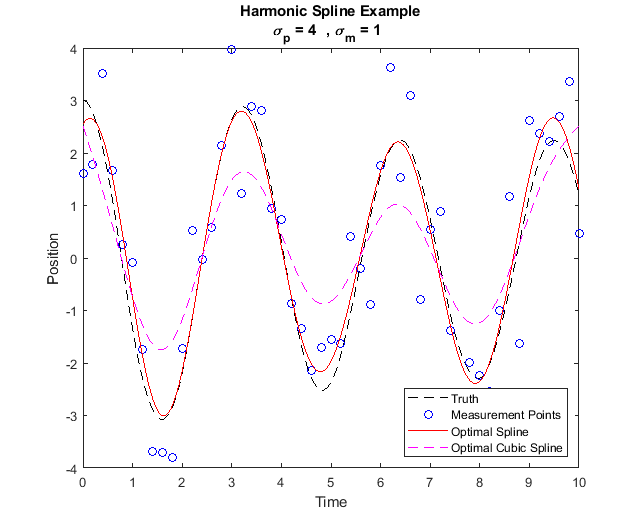}
		\caption{Simulated stochastic harmonic oscillator trajectory with both 
			harmonic and cubic optimal splines for comparison.}
		\label{fig:harm-cubic-comparison}
	\end{figure}
	We have emphasized the adaptive nature of the splines in the new framework. 
	It is interesting to make a direct comparison of the matched spline to that 
	of 
	the cubic spline arising from the simple particle dynamics.
	In Figure \ref{fig:harm-cubic-comparison} we repeat the simulation of the 
	stochastic harmonic oscillator and fit the data using an optimal cubic 
	spline 
	and an optimal modified harmonic spline.
	
	In both examples covered thus far we 
	observed that the spline was linearly dependent on the measurements.
	In the next section we analyze an example with general linear dynamics and 
	measurements subject to Gaussian additive noise.
	These are the class of systems for which the Kalman filter was derived and 
	is 
	known to be optimal.
	We include this illustrative example to
	provide new perspectives on smoothing.
	In particular, we demonstrate that our framework provides an alternative 
	means of smoothing data that circumvents the cumbersome matrix-valued and 
	nonlinear  
	(differential Riccati) equations that arise in Kalman smoothing. 
	Indeed, for linear systems Theorem \ref{thm:Optimization-Conditions} yields 
	linear differential and linear algebraic equations which, once solved, 
	result in
	linear dependence of the splines on measurement data.
	
	\subsection{General Linear Systems}
	\label{subsec:general-kalman-systems}
	We consider a linear system
	\begin{equation}
		\label{eqn:kalman-state}
		\dot{x}(t) = A x(t) + B v(t),
	\end{equation}
	\begin{equation}
		\label{eqn:kalman-measurement}
		y(t) = C x(t) + D w(t),
	\end{equation}
	and Gaussian random processes
	\begin{align}
		\label{eqn:kalman-dyn-dist}
		\rho_v(t,v(t)) = \frac{1}{\sqrt{2 \pi}^{n_v} 
			\mathrm{det}(Q)^{\frac{1}{2}}} 
		e^{-\frac{1}{2} v(t)^T Q^{-1} v(t)}\\
		\label{eqn:kalman-meas-dist}
		\rho_w(t,w(t)) = \frac{1}{\sqrt{2 \pi}^{n_w} 
			\mathrm{det}(R)^{\frac{1}{2}}} 
		e^{-\frac{1}{2} w(t)^T R^{-1} w(t)}.
	\end{align}
	In (\ref{eqn:kalman-dyn-dist}),
	we take $Q$ to be the covariance of $v(t)$.
	Similarly in (\ref{eqn:kalman-meas-dist}) $R$ denotes the covariance of 
	$w(t)$.
	The matrix $A$ is of dimension $n_x \times n_x$, $B$ is of dimension $n_x 
	\times n_v$, $C$ is of dimension $n_y \times 
	n_x$, and 
	$D$ is of dimension $n_y \times n_w$.
	
	As we did in the previous two examples, we apply 
	Theorem\,\ref{thm:Optimization-Conditions} and consider the differential 
	equations that 
	govern the spline,
	\begin{align*}
		&- Q^{-1} v(t) + B^T \lambda(t) = 0\\
		& \dot{\lambda}(t) + A^T \lambda(t) = 0\\
		& \dot{x}(t) - A x(t) - Bv(t) = 0.
	\end{align*}
	This system is solved on $(t_k,t_{k+1})$ by
	\[
	v(t) = QB^T e^{-A^T (t-t_k)} c_{k}^{(\lambda)},
	\]
	\[
	\lambda(t) = e^{-A^T (t-t_k)} c_{k}^{(\lambda)},
	\]
	and
	\[
	x(t) = e^{A(t-t_k)} \big[ \int_{0}^{t-t_k} e^{-A s} B Q B^T e^{-A^T s} ds 
	\hspace{0.1cm} 
	c_{k}^{(\lambda)} 
	+ c_{k}^{(x)}\big],
	\]
	where $c_{k}^{(\lambda)}$ and $c_{k}^{(x)}$ are constant vectors.
	This solution describes the structure of the optimal spline
	for a general linear system with Gaussian processes.
	
	Constraining the constant coefficients of the spline
	requires use of the algebraic equations as we have shown in
	previous examples, and these are written as
	\begin{align*}
		& -Q^{-1} v(t_k) + B^T \lambda(t_k) = 0\\
		& -R^{-1} w(t_k) + D^T \eta(t_k) = 0\\
		& \lambda(t_k) = 0\\
		& A^T \lambda(t_k) + C^T \eta(t_k) + f_0 \lambda(t_k^+) - f_0 
		\lambda(t_k^-) = 0\\
		& y(t_k) - C x(t_k) - D w(t_k) = 0\\
		& \dot{x}(t_k) - A x(t_k) - B v(t_k) = 0,
	\end{align*}
	for each $t_k \in \mathscr{T}_M$.
	Therefore construction of the optimal spline requires only the solution of 
	a 
	\emph{linear} system in this treatment.
	This is not only numerically appealing but we reiterate that it results in 
	closed-form solutions that provide infinite time resolution through a 
	representation as a finite set of constants.
	
	\subsection{An Example with Non-Gaussian Process Noise}
	\label{subsec:non-gaussian-example}
	We move away from Kalman-type examples to further demonstrate benefits 
	of the new framework.
	We revisit the simple particle dynamics 
	under forcing by a non-Gaussian stochastic process.
	A linear filter is non-optimal in this scenario, and indeed the new 
	framework 
	produces a non-linear optimal solution automatically.
	The modified simple particle system is described by
	\[
	\dot{x}(t) = 
	\begin{pmatrix}
		0 & 1\\
		0 & 0
	\end{pmatrix} x(t) + 
	\begin{pmatrix}
		0\\1
	\end{pmatrix}v(t),
	\]
	\[
	y(t) = \begin{pmatrix}
		1 & 0
	\end{pmatrix} x(t) + w(t),
	\]
	\[
	\rho_v(t,v(t)) = c_v 
	e^{-\frac{1}{2}(\frac{v(t)}{\sigma_p})^{2\alpha}},
	\]
	and
	\[
	\rho_w(t,w(t)) = \frac{1}{\sqrt{2 \pi} \sigma_m} 
	e^{-\frac{1}{2 \sigma_m^2}w(t)^2}.
	\]
	The modification introduces a parameter $\alpha$ in the 
	exponent of $\rho_v$, 
	and it is assumed to be a positive integer greater than $1$.
	The leading constant $c_\alpha$ normalizes the distribution and is a 
	function of the parameter, but its specific value is immaterial in our 
	discussion.
	
	The differential equations arising from Theorem 
	\ref{thm:Optimization-Conditions} are written as
	\begin{align}
		\label{eqn:non-linear-process-diff-eq}
		&-\frac{\alpha}{\sigma_p^{2\alpha}} v(t)^{2 \alpha - 1} + 
		\begin{pmatrix}
			0 & 1
		\end{pmatrix} \lambda(t) = 0\\
		\label{eqn:non-linear-lam-de}
		& \dot{\lambda}(t) + \begin{pmatrix}
			0 & 0\\
			1 & 0\\
		\end{pmatrix} \lambda(t) = 0\\
		\label{eqn:non-linear-state-diff-eq}
		& \dot{x}(t) - \begin{pmatrix}
			0 & 1\\
			0 & 0
		\end{pmatrix} x(t) - \begin{pmatrix}
			0\\1
		\end{pmatrix}v(t) = 0.
	\end{align}
	When $\alpha = 1$ we recover the Gaussian result. 
	For $\alpha > 1$ the system is non-linear but admits a closed form solution.
	Equation (\ref{eqn:non-linear-lam-de}) is solved by
	\[
	\lambda(t) = \begin{pmatrix}
		a\\ at + b
	\end{pmatrix},
	\]
	where $a$ and $b$ are arbitrary constants.
	Substitution of $\lambda(t)$ in (\ref{eqn:non-linear-process-diff-eq}) 
	yields
	\[
	v(t) = (\frac{\sigma_p^{2 \alpha}}{\alpha})^{\frac{1}{2 \alpha - 1}}(a t + 
	b)^{\frac{1}{2 \alpha - 1}}.
	\]
	This expression is substituted in (\ref{eqn:non-linear-state-diff-eq}), 
	yielding
	\[
	x(t) = \begin{pmatrix}
		\frac{(2 \alpha - 1)^2}{(4 \alpha - 1)(2 \alpha a^2)} 
		(\frac{\sigma_p^{2 \alpha}}{\alpha})^{\frac{1}{2 \alpha - 1}} (a t + 
		b)^{\frac{4 \alpha - 1}{2 \alpha - 1}} + \frac{c}{a} t + d\\
		\frac{2 \alpha - 1}{2 \alpha a}(\frac{\sigma_p^{2 
				\alpha}}{\alpha})^{\frac{1}{2 \alpha - 1}} (a t + b)^{\frac{2 
				\alpha}{2 
				\alpha -1}} + c
	\end{pmatrix},
	\]
	where $c$ and $d$ are additional arbitrary constants.
	
	This solution determines the spline structure on each of the 
	measurement-free 
	intervals, and therefore we again solve for $4K$ unknown constants.
	These are computed with the algebraic equations from Theorem 
	\ref{thm:Optimization-Conditions}, written as
	\begin{align*}
		& -\frac{\alpha}{\sigma_p^{2\alpha}} v(t_k)^{2 \alpha - 1} + 
		\begin{pmatrix}
			0 & 1
		\end{pmatrix} \lambda(t_k) = 0\\
		& -\frac{1}{\sigma_m^2} w(t_k) + \eta(t_k) = 0\\
		&\lambda(t_k) = 0\\
		&\begin{pmatrix}
			0 & 0\\1 & 0
		\end{pmatrix}\lambda(t_k) + 
		\begin{pmatrix}
			1\\0 
		\end{pmatrix}\eta(t_k) + f_0 \lambda(t_k^+) - f_0 \lambda(t_k^-) = 0\\
		&y(t_k) - \begin{pmatrix}
			1 & 0
		\end{pmatrix} x(t_k) - w(t_k) = 0\\
		&\dot{x}(t_k) - \begin{pmatrix}
			0 & 1\\0 & 0
		\end{pmatrix}x(t_k) - \begin{pmatrix}
			0\\1
		\end{pmatrix}v(t_k) = 0,
	\end{align*}
	for each $t_k \in \mathscr{T}_M$.
	For $\alpha > 1$ the non-linear relationship between the piecewise spline 
	functions and the 
	unknown constants induces a non-linear relationship between the measured 
	data 
	and the optimal estimate when constructing the spline.
	
	\subsection{An Example with Non-linear Dynamics}
	\label{subsec:non-linear-example}
	We conclude the examples with a brief examination of a non-linear dynamical 
	system.
	We consider the non-linearized dynamics of a simple pendulum forced by a 
	Gaussian process. 
	This system is described by
	\[
	\ddot{\theta}(t) + \sin \theta = v(t),
	\]
	where $\theta$ is the angle between the pendulum and direction of the 
	gravitational force.
	We rewrite the system in standard form as
	\begin{align*}
		&\dot{x} = \begin{pmatrix}
			x_2\\ -\sin(x_1)
		\end{pmatrix} + \begin{pmatrix}
			0\\1
		\end{pmatrix}v\\
		&y = \begin{pmatrix}
			1 & 0
		\end{pmatrix} x + w\\
		&\rho_v(t,v(t)) = \frac{1}{\sqrt{2 \pi} \sigma_p} 
		e^{-\frac{1}{2 \sigma_p^2}v(t)^2}\\
		&\rho_w(t,w(t)) = \frac{1}{\sqrt{2 \pi} \sigma_m} 
		e^{-\frac{1}{2 \sigma_m^2}w(t)^2},
	\end{align*}
	where
	\[
	x = \begin{pmatrix}
		\theta\\ \dot{\theta}
	\end{pmatrix}
	\]
	and $y(t)$ models a noisy measurement of $\theta(t)$.
	
	We will compute the optimal spline equations induced by this system, but we 
	do 
	not solve them in the present work as we are not aware of a closed form 
	solution.
	However, we emphasize that our spline equations facilitate the
	application of numerical methods to compute the spline.
	In short, we demonstrate that our framework reduces the enrichment problem 
	to 
	one which simply requires the solution of a system of non-linear 
	differential 
	equations with boundary constraints provided by the algebraic equations.
	
	We note that
	\[
	f(x) = \begin{pmatrix}
		x_2 \\ - \sin (x_1)
	\end{pmatrix},
	\]
	and
	\[
	\frac{\partial f}{\partial x} = \begin{pmatrix}
		0 & 1\\ - \cos (x_1) & 0
	\end{pmatrix},
	\]
	and we use these expressions to simplify notation in the 
	optimization conditions.
	The differential equations provided by Theorem 
	\ref{thm:Optimization-Conditions} are written as
	\begin{align*}
		& -\frac{1}{\sigma_p^2}	v(t) + \begin{pmatrix}
			0 & 1
		\end{pmatrix} \lambda(t) = 0\\
		& \dot{\lambda}(t) + \frac{\partial f}{\partial x}^T \lambda(t) = 0\\
		& \dot{x}(t) - f(x(t)) - \begin{pmatrix}
			0\\ 1
		\end{pmatrix}v(t) = 0.
	\end{align*}
	This is a first order system of non-linear differential equations which 
	describes the general form of the optimal spline.
	The algebraic equations constraining the general solution are written as
	\begin{align*}
		& -\frac{1}{\sigma_p^2}v(t_k) + \begin{pmatrix}
			0 & 1
		\end{pmatrix} \lambda(t_k) = 0\\
		& -\frac{1}{\sigma_m^2}w(t_k) + \eta(t_k) = 0\\
		& \lambda(t_k) = 0\\
		& \frac{\partial f}{\partial x}^T \lambda(t_k) + \begin{pmatrix}
			1\\ 0
		\end{pmatrix} \eta(t_k) + f_0 \lambda(t_k^+) - f_0 \lambda(t_k^-) = 0\\
		& y(t_k) - \begin{pmatrix}
			1 & 0
		\end{pmatrix} x(t_k) - w(t_k) = 0\\
		& \dot{x}(t_k) - f(x(t_k)) - \begin{pmatrix}
			0\\1
		\end{pmatrix} v(t_k) = 0
	\end{align*}
	for each $t_k \in \mathscr{T}_M$.
	At this point the system is fully defined, and one would approach computing 
	its 
	solution using an appropriate numerical technique.
	
	\section{Conclusions and Future Work}
	\label{sec:conclusions}
	In this work we have developed a technique for enriching data using 
	dynamical 
	system and measurement models under additive forcing from stochastic 
	processes.
	We developed an optimization framework that allows the robust modeling 
	of the most general dynamical systems and measurements; one that is not 
	limited 
	to linear dynamics or Gaussian stochastic processes.
	When restricted to the linear Gaussian case our framework naturally, and 
	without prior assumption, renders a linear mapping between the measurements 
	and 
	the optimal spline.
	This differs from Kalman filtering/smoothing theory, where the linearity of 
	this mapping is assumed from the outset.
	
	The capacity to consider more general classes of dynamical systems and 
	stochastic processes when approaching data enrichment in appealing.
	It allows for the consideration of a much larger class of systems which may 
	be encountered in real applications without requiring Gaussian or linear 
	approximations.
	The resulting system of equations in Theorem 
	\ref{thm:Optimization-Conditions} 
	that governs the optimal enrichment in these circumstances will be 
	non-linear 
	in general, and will produce non-linear data filters.
	Investigation and development of new non-linear data filters is an 
	interesting 
	topic of future work.
	Moreover, artificial intelligence (AI) and machine learning (ML) continue 
	to 
	advance rapidly and the need for improved data quality and representations 
	has 
	increased with it.
	The representation of the optimal splines is intrinsically finite and this 
	may 
	have important implications for feature development to support the training 
	of 
	models that will be built using enriched data through modern AI/ML 
	techniques.
	
	The focus of the present work has been on ordinary temporal stochastic 
	dynamical systems, but we are actively pursuing extending the treatment to 
	spatio-temporal systems and other more general problems.
	We are deeply interested in applying this treatment to systems governed by 
	partial, as opposed to ordinary, differential equations.
	Our initial efforts in achieving this goal have begun to bear fruit and we 
	anticipate that future work will produce high performance techniques for 
	performing data enrichment through multiparameter hyper-surfaces as a 
	generalization of the present one-parameter splines.

	\section{Appendix}
	\label{sec:Appendix}
	\subsection{Extending Distributions}
	\label{subsec:extend-instant-dynamics}
	In Section \ref{sec:model-overview-and-prelims} we introduced an infinite 
	product as a route toward defining the distribution extension, but this 
	approach is intractable.
	Instead, we prescribe mathematically meaningful properties for $\mu$, and 
	prove that these induce the expression (\ref{eqn:mu-computation}).
	
	We require three properties: Constancy, Monotonicity, and Geometric 
	Averaging:
	\begin{definition}
		(Constancy) If $\rho = c$ is constant on the interval 
		$\tau$, then $\mu(\tau,\rho) = c$.
	\end{definition}
	\begin{definition}
		(Monotonicity) If $\rho \geq \hat{\rho}$ on the interval $\tau$, then 
		$\mu(\tau,\rho) \geq \mu(\tau,\hat{\rho})$.
	\end{definition}
	\begin{definition}
		(Geometric Averaging) If $\tau_1$ and $\tau_2$ are a partition of the 
		continuous interval $\tau$ such that $\tau = \tau_1\cup \tau_2$, then
		\[
		\mu(\tau,\rho)^{|\tau|} = \mu(\tau_1,\rho)^{|\tau_1|} 
		\mu(\tau_2,\rho)^{|\tau_2|}.
		\]
	\end{definition}
	We show that these properties induce the extension expression used in this 
	work,
	\begin{theorem}
		\label{thm:extension-theorem}
		Assume $\mathrm{ln} \hspace{0.05cm}\rho$ is Lebesgue integrable on 
		$\tau$.
		$\mu(\tau,\rho)$ satisfies constancy, monotonicity, and geometric 
		averaging if and only if
		\[
		\mu(\tau,\rho) = e^{\frac{1}{|\tau|} \int_{\tau} \mathrm{ln} \rho 
			\hspace{0.10cm} dt},
		\]
	\end{theorem}
	and we call Theorem \ref{thm:extension-theorem} the \emph{extension 
		theorem}.
	
	The proof of the theorem uses concepts from Lebesgue measure and 
	Lebesgue integration.
	We do not cover these notions and instead refer to \citep{Wheeden2015} for 
	details.
	
	\begin{proof}
		\underline{\textbf{Forward Direction:}} 
		\newline
		Assume that $\mu(\tau,\rho)$ satisfies constancy, monotonicity, 
		and geometric averaging.
		Let
		\[
		s = \sum_{n=0}^{N-1} s_n \mathrm{1}_{\tau_n}
		\]
		be a simple function on $\tau$, where $\{s_n\}_{n=0}^{N-1}$ is a set of 
		real 
		scalars and
		\[
		\mathrm{1}_{\tau_n} = \begin{cases}
			1 & t \in \tau_n\\
			0 & t \not \in \tau_n.
		\end{cases}
		\]
		The set $\{ \tau_n \}_{n=0}^{N-1}$ is a covering of $\tau$.
		We use the geometric averaging and constancy properties to compute 
		$\mu(\tau,s)$ by writing
		\[
		\mu(\tau,s)^{|\tau|} = \prod_{n=0}^{N-1} s_n^{|\tau_n|},
		\]
		noting that $\mu(\tau_n,s) = s_n$ by construction.
		Taking the logarithm of this expression yields
		\[
		|\tau| \mathrm{ln} \mu(\tau,s) = \sum_{n=0}^{N-1} |\tau_n| 
		\mathrm{ln}(s_n)
		\]
		in general for simple functions.
		The expression on the right hand side of this equation is aligned with 
		the 
		Lebesgue definition of the integral of the logarithm of the simple 
		function,
		\[
		\sum_{n=0}^{N-1} |\tau_n| \mathrm{ln}(s_n) = \int_{\tau} \mathrm{ln} s 
		\hspace{0.10cm} dt.
		\]
		
		The strategy will be to approximate $\mathrm{ln} \rho$ above and below 
		by 
		simple functions which can be used to approximate the integral in the 
		Lebesgue 
		sense to arbitrary precision.
		Given $\epsilon > 0$, since $\mathrm{ln} \rho$ is measurable, there is 
		a simple 
		function $s^* \geq \rho$ 
		such that
		\[
		|\tau| \mathrm{ln} 
		\mu(\tau, \rho) \leq |\tau| \mathrm{ln} \mu(\tau,s^*) = \int_{\tau} 
		\mathrm{ln}s^* dt < 
		\int_{\tau} 
		\mathrm{ln} \rho \hspace{0.10cm} dt + \epsilon.
		\]
		Similarly, there is a simple function $s_* \leq \rho$ such that
		\[
		\int_{\tau} \mathrm{ln} \rho \hspace{0.10cm} dt - \epsilon < 
		\int_{\tau} 
		\mathrm{ln} s_* dt = |\tau| \mathrm{ln} \mu(\tau,s_*) \leq |\tau| 
		\mathrm{ln} \mu(\tau, \rho).
		\]
		Therefore,
		\[
		\Big| |\tau| \mathrm{ln} \mu(\tau,\rho) - \int_{\tau} \mathrm{ln} \rho 
		\hspace{0.10cm} dt \Big| < \epsilon
		\]
		and we conclude that
		\[
		\mathrm{ln} \mu(\tau, \rho) = \frac{1}{|\tau|} \int_{\tau} \mathrm{ln} 
		\rho 
		\hspace{0.10cm}
		dt,
		\]
		as $\epsilon$ is arbitrary.
		
		\noindent \underline{\textbf{Reverse Direction:}}
		\newline
		Assume
		\[
		\mu(\tau,\rho) = e^{\frac{1}{|\tau|}\int_{\tau} \mathrm{ln} \rho 
			\hspace{0.10cm} dt}.
		\]
		We show the properties hold directly.
		For constancy, assume that $\rho = c$ on $\tau$ and consider
		\begin{align*}
			\mu(\tau,\rho) & = e^{\frac{1}{|\tau|}\int_{\tau} \mathrm{ln} \rho 
				\hspace{0.10cm} dt}\\
			& = e^{\frac{1}{|\tau|}\int_{\tau} \mathrm{ln} c 
				\hspace{0.10cm} dt}\\
			& = e^{\frac{\mathrm{ln} c}{|\tau|}\int_{\tau} 
				dt}\\
			& = e^{\mathrm{ln} c}\\
			& = c.
		\end{align*}
		Monotonicity follows from properties of the integral, logarithm, and 
		exponential:
		\[
		\frac{1}{|\tau|}\int_{\tau} \mathrm{ln} \rho dt \geq 
		\frac{1}{|\tau|}\int_{\tau} \mathrm{ln} \hat{\rho} dt
		\]
		for all $\rho \geq \hat{\rho}$.
		Therefore,
		\[
		e^{\frac{1}{|\tau|}\int_{\tau} \mathrm{ln} \rho dt} \geq 
		e^{\frac{1}{|\tau|}\int_{\tau} \mathrm{ln} \hat{\rho} dt}.
		\]
		Finally, we show geometric averaging holds by taking $\tau = \tau_1 
		\cup 
		\tau_2$ for some $\tau_1$ and $\tau_2$, and considering
		\begin{align*}
			\mu(\tau,\rho)^{|\tau|} & = e^{\int_{\tau} \mathrm{ln} \rho 
				\hspace{0.10cm} 
				dt}\\
			& = e^{\int_{\tau_1} \mathrm{ln} \rho 
				\hspace{0.10cm} dt + \int_{\tau_2} \mathrm{ln} \rho 
				\hspace{0.10cm} dt}\\
			& = e^{\int_{\tau_1} \mathrm{ln} \rho 
				\hspace{0.10cm} dt}e^{\int_{\tau_2} \mathrm{ln} 
				\rho \hspace{0.10cm} dt}\\
			& = \mu(\tau_1,\rho)^{|\tau_1|} 
			\mu(\tau_2,\rho)^{|\tau_2|}.
		\end{align*}
		This completes the proof.
	\end{proof}

	\subsection{Application of the Calculus of Variations}
	\label{subsec:Application of the Calculus of Variations}
	We provide a summary of the analysis of problem 
	(\ref{eqn:general-lagrangian-opt-problem}) 
	through demonstration of the calculus of variations.
	For each of the variables in the optimization we introduce a variation. 
	For instance for $x$,
	\[
	\tilde{x} = x + \delta x,
	\]
	where $\tilde{x}$ is the perturbed state, $x$ is the unperturbed state, and 
	$\delta x$ is the variation.
	
	We describe the procedure for a simple, but representative, objective that 
	is 
	written as
	\[
	J(u) = \sum_{k=0}^{K-1} \int_{t_k}^{t_{k+1}} F(u,\dot{u},t) dt + 
	\sum_{k=0}^K G(u(t_k),\dot{u}(t_k),t_k).
	\]
	For the purposes of the summary we use $u$ as a general variable, and in 
	(\ref{eqn:general-lagrangian-opt-problem}) the general variable is 
	constructed 
	by taking
	\[
	u = \begin{pmatrix}
		x\\v\\w\\ \lambda \\ \eta
	\end{pmatrix}.
	\]
	Using standard techniques in the CoV allows us to compute the variation of 
	the 
	integral term
	\[
	\int_{t_k}^{t_{k+1}} F dt
	\]
	via
	\[
	F(u + \delta u,\dot{u} + \delta \dot{u},t)
	= F(u,\dot{u}) + \delta F(u,\dot{u}) + \mathrm{h.o.t},
	\]
	where
	\[
	\delta F(u,\dot{u}) = \frac{\partial F}{\partial u}(u,\dot{u}) \delta u 
	+ \frac{\partial F}{\partial \dot{u}}(u,\dot{u}) \delta \dot{u},
	\]
	to obtain
	\[
	\frac{\partial F}{\partial \dot{u}} \delta 
	u \mid_{t_{k+1}^-} - \frac{\partial F}{\partial 
		\dot{u}} \delta 
	u\mid_{t_{k}^+} - \int_{t_k}^{t_{k+1}} (\frac{d}{dt}\frac{\partial 
		F}{\partial \dot{u}} - \frac{\partial F}{\partial u}) 
	\delta u \hspace{0.10cm} dt.
	\]
	In this expression we use the notation
	\[
	\frac{\partial F}{\partial \dot{u}} \delta 
	u \mid_{t_k^+}	
	= \lim_{t \rightarrow t_k^+} \frac{\partial F}{\partial 
		\dot{u}}(t) \delta 
	u(t)
	\]
	and employ integration-by-parts to remove all temporal derivatives
	from the variation $\delta u$.
	
	The variation of the discrete term
	\[
	G(u(t_k),\dot{u}(t_k),t_k)
	\]
	is written as
	\[
	\frac{\partial G}{\partial u} \delta u \mid_{t_{k}} + 
	\frac{\partial G}{\partial \dot{u}} \delta \dot{u}\mid_{t_{k}}.
	\]
	Taking the sum of these variations (over $k$) and collecting similar terms 
	yields 
	the conditions
	\begin{equation}
		\label{eqn:general-variation-interval}
		\frac{d}{dt}\frac{\partial 
			F}{\partial \dot{u}}(t) - \frac{\partial 
			F}{\partial u}(t) = 0
	\end{equation}
	for all $t \in \mathscr{T} - \mathscr{T}_M$,
	\begin{equation}
		\label{eqn:general-variation-u}
		- \frac{\partial F}{\partial \dot{u}}(t_k^+)
		+ \frac{\partial F}{\partial \dot{u}}(t_k^-)
		+ \frac{\partial G}{\partial u}(t_k) = 0,
	\end{equation}
	and
	\begin{equation}
		\label{eqn:general-variation-u-dot}
		\frac{\partial G}{\partial \dot{u}}(t_k) = 0
	\end{equation}
	for all $t_k \in \mathscr{T}_M$.
	In these expressions we have enforced that the variations must vanish
	for all $\delta u \neq 0$, for 
	otherwise the objective would not be at a critical point and the solution 
	would 
	not be optimal.
	Equation (\ref{eqn:general-variation-interval}),
	when applied to the first term on the right of 
	(\ref{eqn:general-lagrange-problem}),
	generates 
	(\ref{eqn:gen-opt-cond1})--(\ref{eqn:gen-opt-cond3}) in Theorem 
	\ref{thm:Optimization-Conditions}, and
	Equations (\ref{eqn:general-variation-u}) and 
	(\ref{eqn:general-variation-u-dot}),
	when applied to the second term on the right of 
	(\ref{eqn:general-lagrange-problem}),
	generate 
	(\ref{eqn:gen-opt-cond4})--(\ref{eqn:gen-opt-cond9}).
	
	\bibliographystyle{plainnat.bst}
	\bibliography{mlsr.bib}
\end{document}